\documentclass[a4paper]{amsart}
\usepackage{times,amssymb,booktabs,amsmath,mathrsfs,cite}
\usepackage[pagebackref=true]{hyperref}
\renewcommand*{\backref}[1]{}
\renewcommand*{\backrefalt}[4]%
  {{[\tiny\ifcase #1 Not cited.\relax\or Page~#2.\else Pages #2.\fi]}}

\usepackage[svgnames]{xcolor}

\usepackage{aliascnt}
\def\NewTheorem#1{%
  \newaliascnt{#1}{equation}
  \newtheorem{#1}[#1]{#1}
  \aliascntresetthe{#1}
  \expandafter\def\csname #1autorefname\endcsname{#1}
}

\def\equationautorefname~#1\null{(#1)\null}

\newcounter{main}
\theoremstyle{plain}

\swapnumbers
\numberwithin{equation}{section}

\NewTheorem{Proposition}
\NewTheorem{Theorem}
\NewTheorem{Corollary}
\NewTheorem{Lemma}
\NewTheorem{Definition}
\theoremstyle{remark}
\NewTheorem{Remarks}

\newaliascnt{Example}{equation}
\aliascntresetthe{Example}

  {\refstepcounter{equation}\trivlist
   \item[\hskip\labelsep\theequation.~\textbf{Example}\space]
   \ignorespaces
  }{\unskip\nobreak\hfil%
    \penalty50\hskip2em\hbox{}\nobreak\hfil$\Diamond$%
    \parfillskip=0pt\finalhyphendemerits=0\penalty-100\endtrivlist
}

\def\map#1#2{\,{:}\,#1\!\longrightarrow\!#2}
\def\bijection{\overset{\sim}{\longrightarrow}}

\usepackage{xparse,mathtools}
\DeclarePairedDelimiterX{\aset}[1]{\{}{\}}{\setargs{#1}}
\NewDocumentCommand{\setargs}{>{\SplitArgument{1}{|}}m}{\setargsaux#1}
\NewDocumentCommand{\setargsaux}{mm}
   {\IfNoValueTF{#2}{#1} {#1\,\delimsize|\,\mathopen{}#2}}

\let\gdom=\triangleright

\let\gedom=\trianglerighteq

\def\N{\mathbb N}
\def\Q{\mathbb Q}
\def\Z{\mathbb Z}

\def\Sym{\mathfrak S}
\newcommand\Hn[1][n]{\mathscr{H}_{#1}}
\newcommand\HZn[1][n]{\Hn[#1]^\ZZ}
\newcommand\HKn[1][n]{\Hn[#1]^\KK}
\newcommand\HLn[1][n]{\Hn[#1]^\Lambda}
\newcommand\Hlam[1][\blam]{\Hn^{\gdom#1}}
\def\KK{\mathcal{K}}
\def\AA{\mathcal{A}}
\def\ZZ{\mathcal{Z}}

\newcommand\Rn[1][n]{\mathscr{R}^\Lambda_{#1}}
\newcommand\Rlam[1][\blam]{\mathscr{R}^{\gdom#1}}

\newcommand\Parts[1][n]{\mathcal{P}_{#1}}

\newcommand\bi{\mathbf{i}}
\newcommand\bj{\mathbf{j}}

\def\blam{{\boldsymbol\lambda}}
\def\bmu{{\boldsymbol\mu}}

\def\charge{{\boldsymbol\kappa}}

\def\s{\mathfrak s}
\def\t{\mathfrak t}
\def\u{\mathfrak u}
\def\v{\mathfrak v}

\newcommand\Slam[1][\blam]{\mathbb{S}^{#1}}
\newcommand\dSlam[1][\blam]{\mathbb{S}_{#1}}
\def\tlam{\t^\blam}
\def\ilam{\bi^\blam}

\newcommand\Down[1][n]{{\ooalign{\footnotesize\hfil\raise.20ex\hbox{\footnotesize$#1$}\hfil\crcr\square}}}

\usepackage{etoolbox} 
\newcommand{\DeclareMyOperator}[1]{%
  \expandafter\DeclareMathOperator\csname #1\endcsname{#1}
}
\forcsvlist{\DeclareMyOperator}{Res,res,Shape,Std,Add,Rem,Rep,Hom}
\newcommand\Mod{\!\mathop{\text{-Mod}}}
\DeclareMathOperator{\iRes}{\text{$i$}-Res}
\newcommand\defect{\mathop{\textrm{def}}}

\usepackage{tikz}
\newcommand\YoungDiagram[2][\relax]{
  \begin{tikzpicture}[scale=0.43,draw/.append style={thick,black},baseline=-3mm]
    \ifx\relax#1\relax%
    \else 
    \foreach\box in {#1} {
      \filldraw[blue!30]\box rectangle ++(1,1);
    }
    \fi
    \newcount\tableauRow
    \tableauRow=0
    \foreach \diagramCol in {#2} {
      \draw[thin](1,\the\tableauRow)grid ++(\diagramCol,1);
      \global\advance\tableauRow by -1
    }
  \end{tikzpicture}
}
\newcommand\Bigger[2][7]{\left#2\rule{0mm}{#1truemm}\right.}
\def\TriDiagram(#1|#2|#3){\Bigger(%
  \YoungDiagram{#1}\Bigger|\
  \YoungDiagram{#2}\Bigger|\YoungDiagram{#3}
  \Bigger)
}

\usepackage[hyperpageref]{backref}

\def\thetitle{Restricting Specht modules of cyclotomic Hecke algebras}
\hypersetup{pdftitle={\thetitle},
  pdfauthor={Andrew Mathas},
  colorlinks = true,
  linkcolor =ForestGreen,
  anchorcolor = red,
  citecolor = blue,
  urlcolor = blue
}

\def\realMR#1 (#2)#3.{\href{http://www.ams.org/mathscinet-getitem?mr=#1}{MR#1}.}

\begin{document}


  \title{\thetitle}
  \author{Andrew Mathas}
  \address{School of Mathematics and Statistics F07, University of
  Sydney, NSW 2006, Australia.}
  \email{\href{mailto:andrew.mathas@sydney.edu.au}{andrew.mathas@sydney.edu.au}}
  \thanks{This research was supported, in part, by the Australian
  Research Council}

  \begin{abstract}
    This paper proves that the restriction of a Specht module for a
    (degenerate or non-degenerate) cyclotomic Hecke algebra, or KLR
    algebra, of type~$A$ has a Specht filtration.
  \end{abstract}

  \bibliographystyle{andrew}

  \maketitle

\section{Introduction}

  It is a well-known and classical result in the semisimple representation
  theory of the symmetric group $\Sym_n$, for $n\ge1$, that the restriction of
  a Specht module is the direct sum of Specht modules indexed by certain
  partitions of~$n-1$. In general, only the weaker result holds that the
  restriction of a Specht module has a filtration by smaller Specht modules.
  That is, a filtration such that each subquotient is isomorphic to a Specht
  module.

  All of these results generalise to the cyclotomic Hecke algebras of
  type~$A$, which are certain deformations of the group algebra of the
  complex reflection group $\Sym_n\wr(\Z/\ell\Z)$, for $n,\ell\ge1$. In
  the semisimple case, the restriction of a Specht module for a
  cyclotomic Hecke algebra was described by Ariki and Koike~\cite{AK}
  using seminormal forms. In the non-semisimple case, the
  corresponding restriction theorem exists in the literature
  \cite[Proposition~1.9]{AM:simples}, however, the proof given
  in~\cite{AM:simples}, and repeated in
  \cite[Proposition~6.1]{M:ULect}, is incomplete\footnote{Professor
  Ariki disagrees that there is a gap in the original proof in
  \cite{AM:simples}. He says that the result can be easily
  deduced from \cite[Theorem 13.21(2)]{Ariki:book} using the fact
  that the generator of the Specht module is a simultaneous eigenvector
  for the Jucys-Murphy elements.}.
  Ariki~\cite[Lemma~13.2]{Ariki:book} has proved the weaker statement in
  the Grothendieck group. Using the Murphy basis of the Specht
  module~$S^\blam$ for~$\Hn$ it is easy to construct an explicit
  $\Hn[n-1]$-module filtration of~$S^\blam$ that looks like a Specht
  filtration of~$S^\blam$. Unfortunately, as Goodman\footnote{Personal
  communication. See also \cite{GoomanKilgoreTeff}.}
  points out, it is not obvious that the subquotients in this filtration
  are isomorphic to Specht modules. This note gives a self-contained proof
  that the restriction of a Specht module has a Specht filtration (see
  \autoref{T:MurphyRestriction} for a more precise statement).

  Brundan, Kleshchev and Wang \cite[Theorem~4.11]{BKW:GradedSpecht}
  extended the Specht module restriction theorem to the graded setting,
  where this result is an important ingredient in the proof of the
  Ariki-Brundan-Kleshchev categorification
  theorem~\cite[Theorem~5.14]{BK:GradedDecomp}. The proof of the graded
  restriction theorem given in \cite{BKW:GradedSpecht} uses, in an
  essential way, the construction in \cite{AM:simples}. Moreover, the
  arguments of \cite{BKW:GradedSpecht} only apply in the non-degenerate
  case and they require the ground-ring to be a field. We prove the
  graded restriction theorem for graded Specht modules of the cyclotomic
  KLR algebras of type~$A$ defined over an integral domain, thus
  recovering and extending the results of \cite{BKW:GradedSpecht}.

  The layout of this paper is as follows. Section~2 defines the
  cyclotomic Hecke algebras of type~$A$ and sets up the combinatorial
  framework that is used throughout the paper. Section~3 gives a
  self-contained proof of the known branching rule for the Specht
  modules in the semisimple case using the theory of seminormal forms.
  Section~4 using the result from the semisimple case to show that the
  restriction of a Specht module for a cyclotomic Hecke algebra of
  type~$A$ has a Specht filtration, which we describe explicitly.
  Finally, Section~5 introduces the cyclotomic KLR algebras of type~$A$,
  and their graded Specht modules. We show that the restriction of a
  graded Specht module has a filtration by graded Specht modules, up to
  shift, and proves analogous results for the dual graded Specht
  modules.

  The corresponding results for inducing Specht modules can be found in
  \cite{RyomHansen:GradedTranslation,M:SpechtInduction,HuMathas:GradedInduction}.%

\section{Ariki-Koike algebras}\label{S:AKalgebras}
  This section introduces the cyclotomic Hecke algebras of type~$A$.
  Throughout this note we fix positive integers $\ell$ and $n$ and let
  $\Sym_n$ be the symmetric group of degree~$n$. For $1\le r<n$ let
  $s_r=(r,r+1)\in\Sym_n$. Then $s_1,\dots,s_{n-1}$ are the standard
  Coxeter generators of $\Sym_n$.

  We use a slight variation on Ariki and Koike's original definition of
  the cyclotomic Hecke algebras of type~$A$; compare with
  \cite[Definition~3.1]{AK} and \cite[Definition~2.1]{BM:cyc}.  The
  advantage of the definition below is that it allows us to treat the
  so-called degenerate and non-degenerate cases simultaneously, which
  correspond to $\xi^2=1$ and $\xi^2\ne1$, respectively. Previously
  these cases were treated separately.

  Fix a commutative integral domain $\AA$, with $1$.

  \begin{Definition}[%
    \protect{Hu-Mathas\cite[Definition~2.2]{HuMathas:SeminormalQuiver}}]
    \label{D:HeckeAlgebras}
    The \textbf{cyclotomic Hecke
    algebra of type~$A$}, with \textbf{Hecke parameter} $\xi\in\AA^\times$ and
    \textbf{cyclotomic parameters} $Q_1,\dots,Q_\ell\in\AA$, is the unital
    associative $\AA$-algebra $\Hn=\Hn^\AA=\Hn(\AA,\xi,Q_1,\dots,Q_\ell)$
    with generators $L_1,\dots,L_n$, $T_1,\dots,T_{n-1}$ and relations
  {\setlength{\abovedisplayskip}{2pt}
   \setlength{\belowdisplayskip}{1pt}
  \begin{align*}
    \textstyle\prod_{l=1}^\ell(L_1-Q_l)&=0,  &
    (T_r+\xi^{-1})(T_r-\xi)&=0,                  &
    L_{r+1}&=T_rL_rT_r+T_r,
  \end{align*}
  \begin{align*}
      L_rL_t&=L_tL_r, &
    T_rT_s&=T_sT_r\text{ if }|r-s|>1,\\
    T_sT_{s+1}T_s&=T_{s+1}T_sT_{s+1}, &
    T_rL_t&=L_tT_r\text{ if }t\ne r,r+1,
  \end{align*}
  }
  where $1\le r<n$, $1\le s<n-1$ and $1\le t\le n$.
\end{Definition}

The arguments of Ariki and Koike~\cite[Theorem~3.10]{AK} show that $\Hn$
is free as an $\AA$-module with basis
$\aset{L_1^{a_1}\dots L_n^{a_n}T_w|0\le a_1,\dots,a_n<\ell
         \text{ and } w\in\Sym_n}$.
where $T_w=T_{i_1}\dots T_{i_k}$ if $w=s_{r_1}\dots s_{r_k}\in\Sym_n$ is
a reduced expression (that is, $k$ is minimal).

The Ariki-Koike basis theorem implies that there is a natural embedding
of~$\Hn[n-1]$ in~$\Hn$ and that $\Hn$ is free as an $\Hn[n-1]$-module
of rank $\ell n$.  If $M$ is an $\Hn$-module let
\[\Res M\map{\Hn\Mod}\Hn[n-1]\Mod\]
be the restriction functor from category of $\Hn$-modules to
$\Hn[n-1]$-modules. The functor $\Res$ is exact.

The algebra $\Hn$ has another basis that is better adapted to the
representation theory of~$\Hn$.  Recall that a \textbf{partition} of $n$ is a
weakly decreasing sequence $\lambda=(\lambda_1\ge\lambda_2\ge\dots)$ of
non-negative integers such that $|\lambda|=\sum_i\lambda_i=n$. A
\textbf{multipartition}, or $\ell$-partition, of $n$ is an ordered
$\ell$-tuple $\blam=(\lambda^{(1)},\dots,\lambda^{(\ell)})$ of
partitions such that
$|\blam|=|\lambda^{(1)}|+\dots+|\lambda^{(\ell)}|=n$.

Let $\Parts$ be the set of multipartitions of $n$. If
$\blam,\bmu\in\Parts$ then $\blam$ \textbf{dominates} $\bmu$,
written~$\blam\gedom\bmu$, if
\[
  \sum_{c=1}^{k-1}|\lambda^{(c)}|+\sum_{r=1}^s\lambda^{(k)}_r
         \ge\sum_{c=1}^{k-1}|\mu^{(c)}|+\sum_{r=1}^s\mu^{(k)}_r,
\]
for $1\le k\le \ell$ and for all $s\ge1$. Dominance is a partial order
on $\Parts$.

Identify a multipartition $\blam$ with its \textbf{diagram}, which
is the set of nodes
\[\aset{(l,r,c)|1\le\lambda^{(l)}_r\le c\text{ and }1\le l\le \ell}.\]
More generally a \textbf{node} is any element of
$\{1,\dots,\ell\}\times\N^2$, which we order lexicographically.
We picture the diagram of a multipartition as an array of boxes in
the plane, allowing us to talk of the components, rows and columns
of~$\blam$. For example, if $\blam=(3,2|1^3|2^2)$ then
\[
       \blam = \TriDiagram(3,2|1,1,1|2,2)
\]
has three components and $\lambda^{(2)}=(1^3)$ has three rows and one
column.

A \textbf{removable} node of $\blam$ is any node $\rho\in\blam$ such
that $\blam{-}\rho:=\blam\setminus\{\rho\}$ is (the diagram of) a
multipartition. Write $\bmu\to\blam$ if $\bmu$ is obtained from~$\blam$
by removing a removable node. Similarly, a node $\alpha\notin\blam$ is an
\textbf{addable} node if $\blam\cup\aset{\alpha}$ is a multipartition.

The removable nodes of~$\blam$ are totally ordered by the lexicographic
order. As a consequence, the set of multipartitions
$\aset{\bmu\in\Parts[n-1]|\bmu\to\blam}$ is totally ordered by dominance.

If $\blam\in\Parts$ then a $\blam$-tableau is a function
$\t\map{\blam}\aset{1,2,\dots,n}$. Write $\Shape(\t)=\blam$ if $\t$ is a
$\blam$-tableau. For convenience, identify
$\t=(\t^{(1)},\dots,\t^{(\ell)})$ with a labeling of the diagram $\blam$
by elements of $\aset{1,\dots,n}$ in the obvious way. In this way we talk
of the rows, columns and components of~$\t$.

A \textbf{standard $\blam$-tableau} is a $\blam$-tableau $\t$ such that
the entries in each row of $\t^{(k)}$ increase from left to right and the
entries in each column of~$\t^{(k)}$ increase from top to bottom, for
$1\le k\le\ell$. Let $\Std(\blam)$ be the set of standard
$\blam$-tableaux and set $\Std(\Parts)=\bigcup_{\blam\in\Parts}\Std(\blam)$ and
$\Std^2(\Parts)=\bigcup_{\blam\in\Parts}\Std(\blam)\times\Std(\blam)$.

If $\t$ is a standard tableau and $1\le m\le n$ let $\t_{\downarrow m}$
be the subtableau of~$\t$ that contains the numbers $1,2,\dots,m$.  It
is easy to see that~$\t_{\downarrow m}$ is standard because~$\t$ is
standard. Extend the dominance ordering to~$\Std(\Parts)$ by declaring
that $\s\gedom\t$ if $\Shape(\s_{\downarrow
m})\gedom\Shape(\t_{\downarrow m})$, for $1\le m\le n$. As we will need
it often, set $\t_{\downarrow}=\t_{\downarrow(n-1)}$.

Using the observations in the last paragraph, it is easy to verify the
following well-known combinatorial fact, which underpins all of
the results in this note.

\begin{Lemma}\label{L:Bijection}
  Suppose that $\blam\in\Parts$. Then there is a bijection
  \[\Std(\blam)
    \bijection\bigsqcup_{\substack{\bmu\in\Parts[n-1]\\\bmu\to\blam}}
            \Std(\bmu)\qquad(\text{disjoint union}),\]
    given by $\t\mapsto\t_{\downarrow}$.
\end{Lemma}

Let $\tlam$ be the unique standard $\blam$-tableau such that
$\tlam\gedom\t$, for all $\t\in\Std(\blam)$. More explicitly, the
numbers $1,\dots,n$ increase from left to right, and then top to bottom,
in each component of $\tlam$ with the numbers in $\t^{\lambda^{(r)}}$
being less than the numbers in $\t^{\lambda^{(s)}}$ whenever $1\le
r<s\le\ell$.  The symmetric group~$\Sym_n$ acts from the right on the tableaux with
entries in $\{1,2,\dots,n\}$ by permuting their entries.  If $\t$ is a
standard $\blam$-tableau let $d(\t)\in\Sym_n$ be the unique permutation
such that $\t=\tlam d(\t)$.

For each multipartition $\blam\in\Parts$ define
$m_\blam=u_\blam x_\blam$, where
\[
    u_\blam=\prod_{s=2}^\ell\prod_{m=1}^{|\lambda^{(1)}|+\dots+|\lambda^{(s-1)}|}
        (L_m-Q_s)
          \quad\text{ and }\quad x_\blam=\sum_{w\in\Sym_\blam}\xi^{\ell(w)}T_w.
\]
(This corrects the definition of $u_\blam$ given in \cite[\S1.5]{Mathas:Singapore}.)

Let $*$ be the unique anti-isomorphism of~$\Hn$ that fixes all of the
generators in \autoref{D:HeckeAlgebras}. In particular, $T_w^*=T_{w^{-1}}$.
Define the \textbf{Murphy basis} element
$m_{\s\t}=T_{d(\s)}^*m_\blam T_{d(\t)}$, for $\s,\t\in\Std(\blam)$.  By
\cite[Theorem 3.26]{DJM:cyc} and \cite[Theorem~6.3]{AMR},
\begin{equation}\label{E:MurphyBasis}
    \aset{m_{\s\t}|\s,\t\in\Std(\blam)\text{ and }\blam\in\Parts}
\end{equation}
is a cellular basis of $\Hn$, in the sense of Graham and
Lehrer~\cite{GL}. Consequently, if $\Hlam$
is the $\AA$-module spanned by
\(
    \aset{m_{\s\t}|\s,\t\in\Std(\bmu)\text{ for some
             }\bmu\in\Parts \text{ with }\bmu\gdom\blam},
\)
then $\Hlam$ is a two-sided ideal of $\Hn$.

The \textbf{Specht module} $S^\blam$ is the right $\Hn$-submodule of
$\Hn/\Hlam$ generated by $m_\blam+\Hlam$. It follows from
the general theory of cellular algebras that $S^\blam$ is free as an
$\AA$-module with basis $\aset{m_\t|\t\in\Std(\blam)}$, where
$m_\t=m_{\tlam\t}+\Hlam$ for $\t\in\Std(\blam)$. We write
$S^\blam=S^\blam_\AA$ when we want to
emphasize that $S^\blam$ is an $\AA$-module.

Let $M$ be an $\Hn$-module. Then $M$ has a \textbf{Specht filtration}
if there exists a filtration
\[
  0=M_0\subset M_1\subset\dots\subset M_{z}=M
\]
and multipartitions $\blam_1,\dots,\blam_z$ such that
$M_r/M_{r-1}\cong S^{\blam_r}$, for $r=1,\dots,z$. The main results of
this paper say that if $\blam\in\Parts$ then restriction of $S^\blam$ to
$\Hn[n-1]$ has a Specht filtration.

\section{The semisimple case}
  The aim of this note is to understand how the Specht modules of the
  cyclotomic Hecke algebras of type~$A$ behave under restriction. We
  start by considering the easiest case when~$\Hn$ is semisimple. In
  this case, the Specht module restriction theorem is classical,
  building from Young's foundational work on the symmetric group
  from~1901~\cite{QSAI} to the results of Ariki and Koike~\cite{AK}. We
  give the proof here both because it is not very difficult and because
  everything that follows builds upon this result.

  First, we need to define content functions for tableaux. If $k\in\Z$
  define the \textbf{$\xi$-quantum integer} to be the scalar
  \[
            [k]_\xi=\begin{cases*}
              \xi+\xi^3+\dots+\xi^{2k-1}, &if $k\ge0$,\\
              -\xi^{2k}[-k]_\xi,   &if $k<0$.
            \end{cases*}
  \]
  If $\t\in\Std(\Parts)$ is a standard tableau, and $1\le k\le n$, then
  the \textbf{content} of~$k$ in~$\t$ is
  \[c_k(\t)=\xi^{2(c-r)}Q_l+[c-r]_\xi\quad\text{where}\quad(l,r,c)=\t^{-1}(k).\]
  The condition $(l,r,c)=\t^{-1}(k)$ means that $k$ is in row~$r$ and
  column~$c$ of~$\t^{(l)}$.

  The main property of the Murphy basis that we need is that
  the elements $L_1,\dots,L_n$ act upon this basis triangularly. More
  precisely, we have the following:

  \begin{Proposition}[\!{\cite[Proposition 3.7]{JM:cyc-Schaper}
      and~\cite[Lemma~6.6]{AMR}}]\label{P:LkMurphyAction}
    Suppose that $\blam\in\Parts$, $\t\in\Std(\blam)$ and $1\le k\le n$.
    Then
    \[m_\t L_k = c_k(\t)m_\t+\sum_{\s\gdom\t}r_\s m_\s,\]
    for some $r_s\in\AA$.
  \end{Proposition}

  The semisimplicity of $\Hn$ is determined by the \textbf{Poincar\'e
  polynomial} of~$\Hn$:
  \[\displaystyle P_{\Hn}=[1]_\xi[2]_\xi\dots[n]_\xi\prod_{1\le k<l\le\ell}
          \prod_{-n<m<n}(\xi^{2m}Q_k+[m]_\xi-Q_l)\in\AA.\]
  For example, if $\xi^2=1$ then $\Hn\cong\AA\Sym_n$ and $P_{\Hn}=\pm n!$.

  For symmetric groups, the next result is well-known. The
  extension of this result to~$\Hn$ when $\xi^2\ne1$ is due to Ariki~\cite{Ariki:ss}. (The order of
  our exposition is misleading, however, because the most natural way to
  prove the equivalence of parts (a) and (b) in
  \autoref{P:SSconditions} is to use \autoref{T:SeminormalBasis} below.)

  \begin{Proposition}[\protect{Ariki~\cite{Ariki:ss}, %
    Ariki, Mathas, Rui~\cite[Theorem 6.11]{AMR}}]\label{P:SSconditions}
    Suppose that $\AA$ is a field. Then the following are equivalent:
    \begin{enumerate}
      \item The algebra $\Hn$ is semisimple.
      \item If $\s,\t\in\Std(\Parts)$ then $\s=\t$ if and only if
      $c_m(\s)=c_m(\t)$ for $1\le m\le n$.
      \item The Poincar\'e polynomial $P_{\Hn}$ is non-zero.
    \end{enumerate}
  \end{Proposition}

  For the remainder of this section we assume that $P_{\Hn}$ is
  invertible in~$\AA$. The assumption that $P_{\Hn}$ is a unit is useful
  precisely because it implies that $c_m(\t)-c_m(\s)$ is invertible in~$\AA$
  whenever $\s,\t\in\Std(\Parts)$ and $c_m(\s)\ne c_m(\t)$,
  for some $m$. This follows because
  $c_m(\t)-c_m(\s)$ divides $P_{\Hn}$, up to a power of~$\xi$.
  Consequently, if $\t\in\Std(\Parts)$ then
  \[F_\t=\prod_{m=1}^n\prod_{\substack{\s\in\Std(\Parts)\\c_m(\s)\ne c_m(\t)}}
         \frac{L_m-c_m(\t)}{c_m(\t)-c_m(\s)}\quad\in\Hn.
  \]
  If $\s,\t\in\Std(\blam)$ define $f_{\s\t}=F_\s m_{\s\t} F_\t$. Then
  $\aset{f_{\s\t}}$ is a \textbf{seminormal basis} of $\Hn$ in the sense
  of \cite[Definition~3.7]{HuMathas:SeminormalQuiver}. The next result
  summaries that main properties of such bases.

  For the symmetric groups the following result goes back to
  Young~\cite{QSAI}. For the Hecke algebras of types~$A$ and~$B$ this
  was essentially proved by Hoefsmit~\cite{Hoefsmit}. The general cyclotomic
  case was established by Ariki and Koike~\cite{AK}. Here we follow
  the exposition of~\cite{HuMathas:SeminormalQuiver,Mathas:Singapore}.

  \begin{Theorem}\label{T:SeminormalBasis}
    Suppose that $P_{\Hn}\in\AA^\times$.  Then
    $\aset{f_{\s\t}|\s,\t\in\Std(\blam)\text{ for }\blam\in\Parts}$ is a
    basis of $\Hn$ and the following hold:
    \begin{enumerate}
      \item If $(\s,\t)\in\Std^2(\Parts)$ and $1\le m\le n$ then
      \[  f_{\s\t}L_m=c_m(\t) f_{\s\t}\quad\text{and}\quad
          L_m f_{\s\t}=c_m(\s) f_{\s\t}.
      \]
      \item If $(\s,\t)\in\Std^2(\Parts)$ and $1\le r<n$ then
      \[f_{\s\t}T_r=\frac{1+(\xi-\xi^{-1})c_{r+1}(\t)}{c_{r+1}(\t)-c_r(\t)}f_{\s\t}
                 +\alpha_r(\t)f_{\s\v},\]
      where $\v=\t(r,r+1)$ and
      \[\alpha_r(\t)=\begin{cases}
              1,&\text{if }\t\gdom\v,\\
          \frac{(1-\xi^{-1}c_r(\t)+\xi c_r(\v))(1+\xi c_r(\t)-\xi^{-1}c_r(\v))}
          {(c_r(\t)-c_r(\v))(c_r(\v)-c_r(\t))},&\text{otherwise}.
          \end{cases}\]
      \item If $\blam\in\Parts$ then  $S^\blam\cong f_{\s\t}\Hn$, for
      any $\s,\t\in\Std(\blam)$.
      \item If $\AA$ is a field then $\Hn$ is a split semisimple algebra and
      $\aset{S^\blam|\blam\in\Parts}$ is a complete
      set of pairwise non-isomorphic irreducible $\Hn$-modules.
    \end{enumerate}
  \end{Theorem}

  Seminormal bases are most commonly given for Specht modules rather
  than for the regular representation. Parts~(a), (c) and~(d) of
  \autoref{T:SeminormalBasis} are contained in
  \cite[Theorem~3.9]{HuMathas:SeminormalQuiver}, although most of
  \autoref{T:SeminormalBasis} can be deduced from results in
  \cite{AK,AMR}. As shown in
  \cite[Proposition~3.11]{HuMathas:SeminormalQuiver}, part~(b) follows
  because, by \autoref{E:MurphyBasis} and \autoref{P:LkMurphyAction},
  if $(\s,\t)\in\Std^2(\Parts)$ then
  \begin{equation}\label{E:fstExpansion}
      f_{\s\t}=F_\s m_{\s\t}F_\t=m_{\s\t}
           +\sum_{(\u,\v)\in\Std^2(\Parts)} r_{\u\v}m_{\u\v},
      \qquad\text{for some } r_{\u\v}\in\AA,
  \end{equation}
  where $r_{\u\v}\ne0$ only if either
  $\Shape(\u)\gdom\Shape(\s)$ or $(\u,\v)\ne(\s,\t)$,
  $\Shape(\u)=\Shape(\s)$ and $\u\gedom\s$ and $\v\gedom\t$.

  By part~(c) of the theorem, $S^\blam\cong f_{\s\t}\Hn$ for any
  $\s,\t\in\Std(\blam)$. The module $f_{\s\t}\Hn$ has basis
  $\aset{f_{\s\v}|\v\in\Std(\blam)}$ and if $\u\in\Std(\blam)$ then
  $f_{\s\t}\Hn\cong f_{\u\t}\Hn$ with an isomorphism being given by
  $f_{\s\t}\mapsto f_{\u\t}$ by \autoref{T:SeminormalBasis}. In view of
  \autoref{E:fstExpansion}, we can identify $S^\blam$ with the right
  ideal $f_{\tlam\t}\Hn$ of $\Hn$, which has basis
  $\aset{f_{\tlam\t}|\t\in\Std(\blam)}$.

  An almost immediate consequence of \autoref{T:SeminormalBasis} is the
  branching rule for the Specht modules in the semisimple case. Recall
  that $\t_\downarrow=\t_{\downarrow(n-1)}$, for $\t\in\Std(\Parts)$.

  \begin{Corollary}[\protect{Ariki-Koike\cite[Corollary~3.12]{AK}}]
    \label{C:SemisimpleRestriction}
    Suppose that $P_{\Hn}\in\AA^\times$ and fix a multipartition $\blam\in\Parts$.
    Then
    \[\Res S^\blam
    \cong\bigoplus_{\substack{\bmu\in\Parts[n-1]\\\bmu\to\blam}} S^\bmu\]
    as $\Hn[n-1]$-modules. An explicit isomorphism is given by
    $f_{\tlam\t}\mapsto f_{\t^\bmu\t_{\downarrow}}$, where $\t\in\Std(\blam)$
    and $\bmu=\Shape(\t_{\downarrow})$.
  \end{Corollary}

  \begin{proof}
    As above, we may assume that $S^\blam$ has basis
    $\aset{f_{\tlam\t}|\t\in\Std(\blam)}$. If $\bmu\in\Parts[n-1]$ and
    $\bmu\to\blam$ let $S^{\bmu\to\blam}$ be the
    vector subspace of~$S^\blam$ with basis
    $\aset{f_{\tlam\t}|\t_{\downarrow}\in\Std(\bmu)}$. By
    \autoref{L:Bijection}, there are vector space isomorphisms
    \[S^\blam\cong\bigoplus_{\substack{\bmu\in\Parts[n-1]\\\bmu\to\blam}}
                S^{\bmu\to\blam}\
             \cong\bigoplus_{\substack{\bmu\in\Parts[n-1]\\\bmu\to\blam}}
             S^\bmu,\]
    where the second isomorphism is determined by
    $f_{\tlam\t}\mapsto f_{\t^\bmu\t_{\downarrow}}$, for $\t\in\Std(\blam)$
    and where $\bmu=\Shape(\t_{\downarrow})$.  By parts (a)--(c) of
    \autoref{T:SeminormalBasis}, this map restricts to an isomorphism
    $S^{\bmu\to\blam}\cong S^\bmu$ of $\Hn[n-1]$-modules whenever
    $\bmu\in\Parts[n-1]$ and $\bmu\to\blam$. We are done.
  \end{proof}

  We have shown that \autoref{C:SemisimpleRestriction} holds whenever
  $P_{\Hn}$ is invertible in~$\AA$. In particular, this result does not
  require the coefficient ring to be a field.

\section{The non-semisimple case}

  This section proves the first main theorem of this note, which is the
  restriction theorem for the Specht modules over an integral domain.
  This generalises \autoref{C:SemisimpleRestriction}, the branching rule for
  the Specht modules in the semisimple case. The main difference is that
  when $P_{\Hn}$ is invertible then $\Res S^\lambda$ is isomorphic to a
  direct sum of ``smaller'' Specht modules whereas, in general, the
  restriction of a Specht module only has a Specht filtration.

  The next result appears in the literature as
  \cite[Proposition~1.9]{AM:simples}, however, the proof given there is
  incomplete.  The first paragraph of the proof of
  \autoref{T:MurphyRestriction},  which is given below, is essentially
  the same as the proof of this result given in~\cite{AM:simples}. The
  remaining paragraphs complete the argument from \cite{AM:simples} by
  showing that the subquotients in the filtration we construct are
  isomorphic to appropriate Specht modules.

  Recall from \autoref{S:AKalgebras} that $>$ is the lexicographic order
  on the set $\aset{1,\dots,\ell}\times\N^2$ of nodes. Consequently, if
  $A$ and $B$ are removable nodes of $\blam\in\Parts$ then $A>B$ only if
  $A$ is in a latter component of~$\blam$ than $B$, or $A$ and $B$ are
  in the same component and $A$ is in a latter row.

  \begin{Theorem}\label{T:MurphyRestriction}
     Suppose that $\blam\in\Parts$ and let $A_1>\dots>A_z$ be the
     removable nodes of~$\blam$, ordered lexicographically. Then
     $\Res S^\blam$ has an $\Hn[n-1]$-module filtration
     \[
     0=S^{0,\blam}\subset S^{1,\blam}\subset
                     \dots\subset S^{z-1,\blam}\subset S^{z,\blam}=\Res S^\blam
     \]
     such that $S^{r,\blam}/S^{r-1,\blam}\cong S^{\blam-A_r}$ for
     $1\le r\le z$.
   \end{Theorem}

  \begin{proof}
    Let $\bmu_r=\blam-A_r$, for $1\le r\le z$. Then
    $\aset{\bmu_1,\dots,\bmu_z}=\aset{\bmu\in\Parts[n-1]|\bmu\to\blam}$
    and $\bmu_1\gdom\dots\gdom\bmu_z$. Recall from
    \autoref{S:AKalgebras} that $\aset{m_\t|\t\in\Std(\blam)}$ is a basis
    of $S^\blam$. Motivated by the proof of
    \autoref{C:SemisimpleRestriction}, for $1\le r\le z$ define
    $S^{r,\blam}$ to be the $\AA$-submodule of~$S^\blam$ with basis
    $\aset{m_\t|\Shape(\t_{\downarrow})\gedom\bmu_r}$ and set
    $S^{0,\blam}=0$.  Then $S^\blam$ has an $\AA$-module filtration
    \[
         0=S^{0,\blam}\subset S^{1,\blam}\subset \dots\subset S^{z-1,\blam}\subset
                     S^{z,\blam}=S^\blam.
    \]
    We claim that this is the required $\Hn[n-1]$-module filtration of~$S^\blam$.
    By \autoref{P:LkMurphyAction}, if $1\le r\le z$ then the  $\AA$-module
    $S^{r,\blam}$ is stable under the action of $L_1,\dots,L_n$. As remarked in
    \cite{AM:simples}, it is possible to show that $S^{r,\blam}$ is an
    $\Hn[n-1]$-module using \cite[Proposition~3.18]{DJM:cyc}
    (in the cases when ~$\xi^2\ne1$), however, we will prove this directly below.

    Fix $r$ with $1\le r\le z$.
    To complete the proof of the theorem it is enough to show that
    $S^{\bmu_r}\cong S^{r,\blam}/S^{r-1,\blam}$ as $\Hn[n-1]$-modules.
    Let $\theta_r\map{S^{r,\blam}/S^{r-1,\blam}}{S^{\bmu_r}}$ be
    the unique $\AA$-linear map such that
    \[\theta_r(m_\t+S^{r-1,\blam})=m_{\t_{\downarrow}}, \quad
    \text{ for $\t\in\Std(\blam)$ with }\t_{\downarrow}\in\Std(\bmu_r).\]
    By \autoref{L:Bijection}, $\theta_r$ is an $\AA$-module isomorphism. We claim
    that $\theta_r$ is an isomorphism of $\Hn[n-1]$-modules.

    Let $\ZZ=\Z[v,v^{-1},q_1,\dots,q_\ell]$, where $v$,
    $q_1,\dots,q_\ell$ are indeterminates over~$\Z$, and consider the Hecke
    algebra $\HZn=\Hn(\ZZ,v,q_1,\dots,q_\ell)$. Then
    $\Hn\cong\HZn\otimes_\ZZ\AA$ and
    $S^\blam_\AA\cong S^\blam_\ZZ\otimes_\ZZ\AA$, where we
    consider $\AA$ as a $\ZZ$-module by letting~$v$ act as
    multiplication by~$\xi$ and $q_k$ as multiplication by $Q_k$, for $1\le k\le\ell$.
    Similarly, let $\HZn[n-1]$ be the obvious subalgebra of $\HZn$ such
    that $\Hn[n-1]\cong\HZn[n-1]\otimes_\ZZ\AA$.
    By definition, the map $\theta_r$ commutes with base change. Therefore, in order
    to show that~$\theta_r$ is an $\Hn[n-1]$-module homomorphism it is
    enough to consider the case when $\AA=\ZZ$. Let
    $\KK=\Q(v,q_1,\dots,q_\ell)$ be the field of fractions of~$\ZZ$.  Then
    $S^\blam_\KK\cong S^\blam_\ZZ\otimes_\ZZ\ZZ$, so by considering
    $S^\blam_\ZZ$ as a $\ZZ$-submodule of $S^\blam_\KK$ it follows that
    $\theta_r$ is a homomorphism of~$\Hn[n-1]^\ZZ$-modules if and only if the
    induced map~$\theta_r\otimes1_\KK$ over~$\KK$ is a homomorphism of
    $\HKn[n-1]$-modules.  Hence, we are reduced to the case when
    $\AA=\KK$.

    Assume now that $\AA=\KK$ and $\Hn=\HKn$. Then $P_{\HKn}\ne0$, so
    $\HKn$ has a seminormal basis $\aset{f_{\s\t}}$ by \autoref{T:SeminormalBasis}.
    For $\t\in\Std(\blam)$ set $f_\t=f_{\tlam\t}+\Hlam$. Then
    $\aset{f_\t|\t\in\Std(\blam)}$ is a basis of~$S^\blam_\KK$. By
    \autoref{E:fstExpansion}, if $\t\in\Std(\blam)$ then
    \[f_\t = m_\t+\sum_{\s\gdom\t}a_\s m_\s \qquad\text{for some }a_\s\in\KK.\]
    Consequently, $S^{r,\blam}_\KK$ has basis
    $\aset{f_\t|\t\in\Std(\blam)\text{ and }\Shape(\t_\downarrow)\gedom\bmu_r}$.
    In particular, $S^{r,\blam}_\KK$ is an $\Hn[n-1]$-submodule of $S^\blam_\KK$
    by \autoref{T:SeminormalBasis}.  Although we do not need this, this
    also implies that $S^{r,\blam}_\KK\cong
    S^{\bmu_1\to\blam}_\KK\oplus\dots\oplus S^{\bmu_r\to\blam}_\KK$ as
    $\Hn[n-1]$-modules, where we use the notation from the proof of
    \autoref{C:SemisimpleRestriction}.

    By \autoref{C:SemisimpleRestriction} the $\KK$-module isomorphism
    determined by
    \[
       \hat\theta\map{S^\blam_\KK}\bigoplus_{\substack{\bmu\in\Parts[n-1]\\\bmu\to\blam}}
             S^\bmu_\KK; f_\t\mapsto f_{\t_{\downarrow}},
    \]
    is an isomorphism of $\Hn[n-1]^\KK$-modules. Moreover, since the
    transition matrix between the Murphy basis $\aset{m_\t}$ and the
    seminormal basis $\aset{f_\t}$ of $S^\blam_\KK$ is unitriangular, it
    follows by downwards induction on dominance that
    $\hat\theta(m_\t)=m_{\t_{\downarrow}}$, for all $\t\in\Std(\blam)$. By the
    last paragraph, $\hat\theta$ induces an $\Hn[n-1]^\KK$-homomorphism
    $\hat\theta_r\map{S^{r,\blam}_\KK/S^{r-1,\blam}_\KK}S^{\bmu_r}_\KK$
    given by
    \[\hat\theta_r(f_\t+S^{r-1,\blam}_\KK) = \hat\theta(f_\t)=f_{\t_{\downarrow}},\qquad
              \text{for $\t\in\Std(\blam)$ with }\Shape(\t_{\downarrow})=\bmu_r.
    \]
    Therefore,
    $\hat\theta_r(m_\t+S^{r-1,\blam}_\KK)=\hat\theta(m_\t)
              =m_{\t_{\downarrow}}=\theta_r(m_\t+S^{r-1,\blam}_\KK)$.
    Hence, $\theta_r=\hat\theta_r$ is an $\HKn[n-1]$-module homomorphism
    as we wanted to show. Moreover, there is an isomorphism of
    $\HKn[n-1]$-modules $S^{\bmu_r}_\KK\cong S^{r,\blam}_\KK/S^{r-1,\blam}_\KK$,
    for $1\le r\le z$, so this completes the proof.
  \end{proof}

  In level one (that is, when $\ell=1$), a more complicated proof of
  \autoref{T:MurphyRestriction}, which builds on a sketch by the author, can be
  found in \cite{GoomanKilgoreTeff}.

\section{Graded Specht modules}

  In this final section we consider the analogue of
  \autoref{T:MurphyRestriction} for the graded Specht modules of the
  corresponding cyclotomic KLR (or quiver) Hecke algebra~$\Rn$. The
  cyclotomic KLR Hecke algebras~\cite{KhovLaud:diagI,Rouq:2KM} are a
  family of $\Z$-graded algebras that arise in the categorification of
  quantum groups. For quivers of type~$A$, Brundan and
  Kleshchev~\cite{BK:GradedDecomp} have shown that over a field the
  cyclotomic quiver Hecke algebras are isomorphic to the cyclotomic
  Hecke algebras.

Let us set the notation. Fix
$e\in\aset{2,3,4,\dots}\cup\aset{\infty}$ and let $\Gamma_e$ be the quiver
with vertex set $I=\Z/e\Z$ (put $e\Z=\aset{0}$ when $e=\infty$), and with edges
$i\to i+1$, for $i\in I$. To the quiver $\Gamma_e$ we attach the
standard Lie theoretic data of a Cartan matrix $(c_{ij})_{i,j\in I}$,
fundamental weights $\aset{\Lambda_i|i\in I}$, positive weights
$P^+=\sum_{i\in I}\N\Lambda_i$ and  positive roots
$Q^+=\bigoplus_{i\in I}\N\alpha_i$. Let $(\cdot,\cdot)$ be the bilinear form satisfying
\[(\alpha_i,\alpha_j)=c_{ij}\qquad\text{and}\qquad
          (\Lambda_i,\alpha_j)=\delta_{ij},\qquad\text{for }i,j\in I.\]
Let $Q^+_n=\aset{\beta\in Q^+|\sum_{i\in I}(\Lambda_i,\beta)=n}$ and set
$I^\beta=\aset{\bi\in I^n|\beta=\alpha_{i_1}+\dots+\alpha_{i_n}}$
for $\beta\in Q^+_n$.

A \textbf{graded $\AA$-module} or, more accurately a $\Z$-graded
$\AA$-module, is an $\AA$-module $M$ with a decomposition
$M=\bigoplus_{d\in\Z}M_d$ such that each summand $M_d$ is $\AA$-free and
of finite rank. We always assume that only finitely many of the $M_d$
are non-zero. A graded $\AA$-algebra is an $\AA$-algebra
$A=\bigoplus_d A_d$ that is a graded module such that $A_cA_d\subseteq
A_{c+d}$, for all $c,d\in\Z$. A graded $A$-module is an $A$-module that
is graded and $M_cA_d\subseteq M_{c+d}$, for all $c,d\in\Z$.  If $M$ is
a graded $A$-module let $\underline{M}$ be the module obtained by
forgetting the grading. If $s\in\Z$ let $q^sM$, where $q$ is an
indeterminate, be the graded $A$-module obtained by shifting the grading
by~$s$, so that $(q^sM)_d=M_{d-s}$. Two graded $A$-modules are
isomorphic if there is a degree preserving isomorphism between them. For
example, $M\cong q^sM$ if and only if $s=0$.

\begin{Definition}\label{D:CycKLR}
  Fix $e\in\aset{2,3,4,\ldots}\cup\{\infty\}$ and suppose that $\Lambda\in P^+$ and
  $\beta\in Q^+$. The \textbf{cyclotomic quiver Hecke algebra}
  $\Rn[\beta]$ is the unital associative $\AA$-algebra  with generators
  \begin{equation*}
    \{\psi_1,\ldots,\psi_{n-1}\}\cup \{y_1,\ldots,y_n\}
         \cup \aset{e(\bi)| \bi\in I^\beta}
  \end{equation*}
  and relations
  {\setlength{\abovedisplayskip}{2pt}
   \setlength{\belowdisplayskip}{1pt}
  \begin{xalignat*}{3}
    y_1^{(\Lambda,\alpha_{i_1})}e(\bi)&= 0,& e(\bi)e(\bj)&= \delta_{\bi\bj}e(\bi),
        &\textstyle\sum_{\bi\in I^\beta}e(\bi)&= 1,\\
  y_re(\bi)&= e(\bi)y_r,& \psi_re(\bi)&= e(s_r\cdot\bi)\psi_r,& y_ry_s&= y_sy_r,
  \end{xalignat*}
  \begin{xalignat*}{2}
    \psi_ry_{r+1}e(\bi)&= (y_r\psi_r+\delta_{i_r i_{r+1}})e(\bi),
  & y_{r+1}\psi_re(\bi)&= (\psi_ry_r+\delta_{i_r i_{r+1}})e(\bi),
  \end{xalignat*}
  \begin{align*}
  \psi_ry_s&= y_s\psi_r,&\text{if }s\neq r,r+1,\\
  \psi_r\psi_s&= \psi_s\psi_r,&\text{if }|r-s|>1,
  \end{align*}
  \begin{align*}
    \psi_r^2e(\bi)&= Q_{i_ri_{r+1}}(y_r,y_{r+1})e(\bi),\\
  (\psi_r\psi_{r+1}\psi_r-\psi_{r+1}\psi_r\psi_{r+1})e(\bi)
  &=\delta_{i_r,i_{r+2}}
     \frac{Q_{i_{r}i_{r+1}}(y_{r+2},y_{r+1})-Q_{i_{r}i_{r+1}}(y_r,y_{r+1})}
        {y_{r+2}-y_r}
  \end{align*}
  }%
  for $\bi,\bj\in I^\beta$ and all admissible $r$ and $s$ and where
  \[
      Q_{ij}(u,v) = \begin{cases*}
           (u-v)(v-u),& if $i\leftrightarrows j$,\\
           (u-v),& if $i\rightarrow j$,\\
           (v-u),& if $i\leftarrow j$,\\
           0,& if $i=j$,\\
           1,&otherwise.
      \end{cases*}
  \]
  For $n\ge0$, define
  $\displaystyle \Rn =\bigoplus_{\beta\in Q^+_n}\Rn[\beta]$.
\end{Definition}

The algebra $\Rn$ is a $\Z$-graded algebra with degree function
determined by
\[ \deg e(\bi)=0,\quad
   \deg y_r=2\quad\text{and}\quad
   \deg\psi_se(\bi)=(\alpha_{i_r},\alpha_{i_{r+1}}),
\]
for $\bi\in I^n$, $1\le r\le n$ and $1\le r<n$.

Inspecting the relations in \autoref{D:CycKLR}, $\Rn$ has a unique
anti-isomorphism $\diamond$ that fixes each of the generators.

To connect these algebras with the cyclotomic Hecke algebras,
let $e\in\aset{2,3,4,\dots}\cup\aset{\infty}$ be minimal such that
$[e]_\xi=0$, or set $e=\infty$ if $[k]_\xi\ne0$ for $k>0$.

For $\Lambda\in P^+$ choose a \textbf{multicharge}
$\charge=(\kappa_1,\dots,\kappa_\ell)\in\Z^\ell$ such that
\[ (\Lambda,\alpha_i) = \#\aset{1\le k\le\ell|\kappa_k\equiv i\pmod e}
                     \qquad\text{for all }i\in I.\]
Let $\HLn$ be the cyclotomic Hecke algebra with parameters
$Q_r=[\kappa_r]_\xi$, for $1\le r\le\ell$. Up to isomorphism, the
algebra $\HLn$ depends only on~$\Lambda$ and not on the choice of
multicharge.

\begin{Theorem}[Brundan and Kleshchev~\cite{BK:GradedKL}]\label{T:KLR}
  Suppose that $\AA$ is a field. Then $\HLn\cong\Rn$.
\end{Theorem}

The isomorphism of \autoref{T:KLR} holds only over a field even though
the cyclotomic KLR algebra $\Rn$ is defined over any ring. For each
$\blam\in\Parts$ there is a graded Specht module $\Slam$ such that
$\underline{\mathbb{S}}^\blam\cong S^\blam$ as (ungraded)
$\HLn$-modules. To define these, and to prove a graded version of
\autoref{T:MurphyRestriction}, we need some notation.

Fix $i\in I$. A node $A=(l,r,c)$ is an $i$-node if $i=\kappa_l+c-r+e\Z$.
Similarly, if $\t\in\Std(\Parts)$ and $1\le k\le n$ then the \textbf{$e$-residue}
of~$k$ in~$\t$ is
\[  \res_k(\t)=\kappa_l+c-r+e\Z\quad\text{if}\quad\t(l,r,c)=k. \]
The \textbf{residue sequence} of $\t$ is
$\res(\t)=(\res_1(\t),\dots,\res_n(\t))\in I^n$. Set
$\ilam=\res(\tlam)$.

For $i\in I$ and $\blam\in\Parts$ let $\Add_i(\blam)$ be the set of
addable $i$-nodes for~$\blam$ and let $\Rem_i(\blam)$ be the set of
removable $i$-nodes. Order $\Add_i(\blam)$ and $\Rem_i(\blam)$ lexicographically
by~$<$.

  \begin{Definition}[\protect{%
    Brundan, Kleshchev and Wang~\cite[\S3.5]{BKW:GradedSpecht}}]
    \label{D:TableauxDegree}
    Suppose that $\blam\in\Parts$ and let~$A$ be an addable or removable $i$-node
    of~$\blam$. Define integers
    \begin{align*}
       d_A^\blam&=\#\aset{B\in\Add_i(\blam)|B>A}-\#\aset{B\in\Rem_i(\blam)|B>A},\\
       d^A_\blam&=\#\aset{B\in\Add_i(\blam)|B<A}-\#\aset{B\in\Rem_i(\blam)|B<A}.
    \end{align*}
    If $\t\in\Std(\blam)$ define the \textbf{degree} of~$\t$ inductively
    by setting $A=\t^{-1}(n)$ and
    \[\deg\t=\begin{cases}
                \deg \t_{\downarrow} +d_A^\blam, &\text{if }n>0,\\
                0, &\text{if }n=0.
              \end{cases}
    \]
  \end{Definition}

  Following~\cite[Definition~4.9]{HuMathas:GradedCellular}, for
  $\blam\in\Parts$ define $d^\blam_k=\deg\tlam_{\downarrow
  k}-\deg\tlam_{\downarrow(k-1)}$, for $k=1,2,\dots,n$. Define
  $e^\blam = e(\ilam)$ and
  $y^\blam = y_1^{d^\blam_1}y_2^{d^\blam_2}\dots y_n^{d^\blam_n}$.
  For $\s,\t\in\Std(\blam)$ set
  \[\psi_{\s\t} = \psi_{d(\s)}^\diamond e^\blam y^\blam\psi_{d(\t)},\]
  where for each $w\in\Sym_n$ we fix an arbitrary reduced expression
  $w=s_{r_1}\dots s_{r_k}$ and define
  $\psi_w=\psi_{r_1}\dots\psi_{r_k}\in\Rn$. In general, the elements
  $\aset{\psi_{\s\t}}$ depend upon the choices of reduced expression for
  the elements of~$\Sym_n$. By definition,
  $(\psi_{\s\t})^\diamond=\psi_{\t\s}$.

  Hu and the author~\cite{HuMathas:GradedCellular} proved that if $\AA$
  is a field then $\aset{\psi_{\s\t}}$ is a graded cellular basis
  of~$\Rn$, where graded cellularity is a natural extension of the
  definition of cellular algebras to the graded setting. In particular,
  the existence of a graded
  cellular basis automatically gives graded Specht modules. In fact,
  $\Rn$ is a graded cellular algebra over any integral domain.

  \begin{Theorem}[\protect{Li~\cite[Theorem~1.1]{GeLi:IntegralKLR}}]
    \label{T:IntegralKLR}
     Suppose that $\Lambda\in P^+$ and let $\AA$ be an integral
     domain. Then
     $\aset{\psi_{\s\t}|\s,\t\in\Std(\blam)\text{ for some }\blam\in\Parts}$
     is a graded cellular basis of~$\Rn$. In particular,
     $\deg\psi_{\s\t}=\deg\s+\deg\t$, for $\s,\t\in\Std(\blam)$ and
     $\blam\in\Parts$.
  \end{Theorem}

  Following~\cite{HuMathas:GradedCellular}, we now define the graded Specht
  modules although we need to be slightly careful to choose the
  ``correct'' grading shift. In more detail, if $\blam\in\Parts$ let
  $\Rlam$ be the submodule of $\Rn$ spanned by
  $\aset{\psi_{\s\t}|\s,\t\in\Std(\bmu)\text{ where }\bmu\gdom\blam}$.
  Then $\Rlam$ is a two-sided ideal of~$\Rn$. The
  \textbf{graded Specht module} is the $\Rn$-module
  \[ \Slam=q^{-\deg\tlam}(\psi_{\tlam\tlam}+\Rlam)\Rn. \]
  Set $\psi_\t=\psi_{\tlam\t}+\Rlam\in\Slam$. Then
  $\aset{\psi_\t|\t\in\Std(\blam)}$ is a homogeneous basis of $\Slam$,
  where the grading is given by $\deg\psi_\t=\deg\t$, for $\t\in\Std(\blam)$.
  By \autoref{T:IntegralKLR}, the Specht module $\Slam$ is defined over
  an integral domain. The integral graded Specht module can also be
  defined using \cite{KMR:UniversalSpecht}, which gives homogeneous
  relations for the graded Specht modules.

  Importantly, if $\AA$ is a field and $\blam\in\Parts$ then
  $\underline{\mathbb{S}}^\blam\cong S^\blam$ as (ungraded)
  $\HLn$-modules, so it is natural to ask for a graded analogue of
  \autoref{T:MurphyRestriction}.  In the non-degenerate case
  ($\xi^2\ne1)$, and when working over a field, this result appears as
  Brundan, Kleshchev and Wang \cite[Theorem~4.11]{BKW:GradedSpecht},
  however, their proof relies on the arguments of~\cite{AM:simples}. We
  give a self-contained proof of this result that applies over an
  integral domain.

  \begin{Theorem}\label{T:GradedRestriction}
     Suppose that $\AA$ is an integral domain, $\blam\in\Parts$ and let
     $A_1>\dots>A_z$ be the removable nodes of~$\blam$,
     ordered lexicographically. Then $\Res\Slam$ has an $\Rn[n-1]$-module
     filtration
     \[
        0=\Slam[0,\blam]\subset \Slam[1,\blam]\subset
                     \dots\subset \Slam[z-1,\blam]\subset \Slam[z,\blam]=\Res\Slam
     \]
     such that
     $\Slam[r,\blam]/\Slam[r-1,\blam]\cong q^{d^\blam_{A_r}}\Slam[\blam{-}A_r]$,
     for $1\le r\le z$.
  \end{Theorem}

  The special case of \autoref{T:GradedRestriction} when~$\AA$ is a field,
  and $e$ is not equal to the characteristic, recovers
  \cite[Theorem~4.11]{BKW:GradedSpecht}.

  \begin{proof}
    The proof follows the same strategy as the proof of
    \autoref{T:MurphyRestriction}. For $1\le r\le n$ let
    $\bmu_r=\blam{-}A_r$ and define $\Slam[r,\blam]$ to be the
    $\AA$-submodule of~$\Slam$ with basis
    \[ \aset{\psi_\t|\Shape(\t_{\downarrow})\gedom\bmu_r}. \]
    Then, as an $\AA$-module, $\Slam$ has a filtration
    $0=\Slam[0,\blam]\subset \Slam[1,\blam]\subset
                        \dots\subset \Slam[z,\blam]=\Slam.$
    By \autoref{L:Bijection} and \autoref{D:TableauxDegree},
    there are graded $\AA$-module isomorphisms
    \[\theta_r\map{\Slam[r,\blam]/\Slam[r-1,\blam]}
           q^{d^\blam_{A_r}}\Slam[\blam{-}A_r];
           \psi_\t+\Slam[r-1,\blam]\mapsto \psi_{\t_{\downarrow}}
    \]
    for $1\le r\le z$. By definition, $\theta_r$ is homogeneous,
    however, it is not clear that it is a~$\Rn[n-1]$-module homomorphism.

    Fix $r$ with $1\le r\le z$. We claim that $\theta_r$ is an
    isomorphism of graded $\Rn[n-1]$-modules. By \autoref{T:IntegralKLR}, or
    \cite[Corollary~6.24]{KMR:UniversalSpecht},
    $\Slam_\AA\cong\Slam_\Z\otimes_\Z\AA$. By definition, $\theta_r$
    commutes with base change so it is enough to prove the theorem in the
    special case when $\AA=\Z$.  Let $\xi$ be a primitive $e$th root of
    unity, or an indeterminate over~$\Q$ if $e=\infty$, and
    let~$\KK=\Q(\xi)$. Then $\Slam_\Z$ embeds in
    $\Slam_\KK\cong\Slam_\Z\otimes_\Z\KK$, so it is enough to prove the
    theorem over~$\KK$. Now, by \autoref{T:KLR}, $\Rn\otimes_\Z\KK\cong\Hn^\KK$,
    so we are reduced to showing that $\theta_r$ is an $\HKn$-module
    homomorphism. To complete the proof we need one more observation: by
    \cite[Proposition~4.5]{BKW:GradedSpecht}, or
    \cite[Lemma~5.4]{HuMathas:GradedCellular}, there exist non-zero
    scalars $c_\t\in\KK$ such that
    \[ \psi_\t = c_\t
    m_\t+\sum_{\substack{\v\in\Std(\blam)\\\v\gdom\t}}r_\v m_\v,
              \qquad\text{ for }r_\v\in\KK.\]
    That is, the transition matrix between the $\psi$-basis of
    $S^\blam_\KK$ and the Murphy basis is triangular. It follows that
    the map~$\theta_r$, defined above, coincides exactly with the map
    $\theta_r$ defined in the proof of \autoref{T:MurphyRestriction}. In
    particular, $\theta_r$ is a (homogeneous) $\HKn$-module homomorphism,
    as we needed to show.
  \end{proof}

  \begin{Remarks}\hspace*{0.8\textwidth}
    \begin{enumerate}
      \item The argument above proves \autoref{T:GradedRestriction} by
      reducing to the situation considered in
      \autoref{T:MurphyRestriction} that, in turn, reduces to the
      semisimple case of \autoref{C:SemisimpleRestriction}. Using the
      framework developed in \cite{HuMathas:SeminormalQuiver} it is
      possible to deduce \autoref{T:GradedRestriction} directly from
      \autoref{C:SemisimpleRestriction}.
      \item In the special case when~$\AA$ is a field,
      \autoref{T:GradedRestriction} can be deduced almost directly from
      \autoref{T:MurphyRestriction} since the transition matrix between
      the $\psi$-basis and the Murphy basis of the Specht module is
      unitriangular. We use the integral form of~$\Slam$ in the proof
      of \autoref{T:GradedRestriction} only to obtain the result over an
      integral domain.
      \item Because of \autoref{D:HeckeAlgebras}, the proof of
      \autoref{T:GradedRestriction} simultaneously treats both the
      degenerate and non-degenerate cases. Only the non-degenerate case
      was considered in \cite{BKW:GradedSpecht}.
    \end{enumerate}
  \end{Remarks}

  Fix $i\in I$. If $\bj\in I^{n-1}$ let $\bj\vee i=(j_1,\dots,j_{n-1},i)\in I^n$
  and define
  \[
         e_{n-1,i}=\sum_{\bj\in I^{n-1}}e(\bj\vee i)\in\Rn.
  \]
  Using the relations, and \autoref{T:KLR}, there is a non-unital homogeneous
  algebra embedding $\Rn[n-1]\hookrightarrow\Rn$ given by
  \[ e(\bj)\mapsto e(\bj\vee i), \quad y_r\mapsto e_{n-1,i}y_r,
       \quad\text{and}\quad \psi_s\mapsto e_{n-1,i}\psi_s,
  \]
  for $\bj\in I^{n-1}$, $1\le r<n$ and $1\le s<n-1$.
  If $i\in I$ and $M$ is an $\Rn$-module define
  $\iRes M = Me_{n-1,i}$. Then there is an isomorphism of functors
  \[
     \Res\cong\bigoplus_{i\in I}\iRes
  \]
  corresponding to the identity $1_{\Rn}=\sum_{i\in I}e_{n-1,i}$.
  By \cite{BK:GradedKL,Brundan:degenCentre,LM:AKblocks}, the blocks of
  $\Rn[n-1]$ are labelled by $Q^+_{n-1}$.  For $\beta\in Q^+_{m}$ let
  $\Parts[\beta]=\aset{\bmu\in\Parts[m]|\bi^\bmu\in I^\beta}$, for $m\ge0$.

  Applying the idempotent $e_{n-1,i}$ to \autoref{T:GradedRestriction}
  gives the following:

  \begin{Corollary}\label{C:GradedRestriction}
     Let $\AA$ be an integral domain. Suppose that $i\in I$ and
     $\blam\in\Parts[\beta+\alpha_i]$, where $\beta\in Q^+_{n-1}$. Let
     $B_1>\dots>B_y$ be the removable $i$-nodes of~$\blam$, ordered
     lexicographically. Then $\iRes\Slam$ has an $\Rn[\beta]$-module
     filtration
     \[     0=\Slam[0,\blam]\subset\Slam[1,\blam]\subset
                    \dots\subset \Slam[y-1,\blam]\subset \Slam[y, \blam]=\iRes\Slam
     \]
    such that
    $\Slam[r,\blam]/\Slam[r-1,\blam]\cong q^{d^\blam_{B_r}}\Slam[\blam{-}B_r]$
    as $\Rn[\beta]$-modules for $1\le r\le y$.
  \end{Corollary}


  Let $\Rep(\Rn)$ be the Grothendieck group of $\Rn$. If $M$ is an
  $\Rn$-module let $[M]$ be its image in $\Rep(\Rn)$. Then $\Rep(\Rn)$ is
  free as a $\Z[q,q^{-1}]$-module, where $[qM]=q[M]$.
  By \autoref{T:GradedRestriction}, if $i\in I$ then
  \[[\iRes S^\blam]=\sum_{r=1}^yq^{d^\blam_{B_r}}[\Slam[\blam-B_r]].\]
  This result is used in ~\cite[Proposition~5.5]{BK:GradedDecomp} to
  show that the cyclotomic Hecke algebras categorify the highest weight
  representations of the quantum affine special linear group.

  Finally, the papers
  \cite{HuMathas:GradedInduction,KMR:UniversalSpecht} also define
  \textbf{dual graded Specht modules} $\dSlam$, for $\blam\in\Parts$.
  These modules can be defined in exactly the same way as the graded
  Specht modules using a ``dual'' graded cellular basis
  $\aset{\psi'_{\s\t}}$ of~$\Rn$.  This allows all of the arguments
  in this section to be repeated, essentially word-for-word, to prove
  analogues of \autoref{T:GradedRestriction} and
  \autoref{C:GradedRestriction} for the dual graded Specht modules. We
  state only the dual Specht module version of
  \autoref{C:GradedRestriction}.

  \begin{Corollary}\label{C:dualGradedRestriction}
     Let $\AA$ be an integral domain. Suppose that $i\in I$ and
     $\blam\in\Parts[\beta+\alpha_i]$, where $\beta\in Q^+_{n-1}$. Let
     $B_y>\dots>B_1$ be the removable $i$-nodes of~$\blam$, ordered
     lexicographically. Then $\iRes\dSlam$ has an $\Rn[\beta]$-module
     filtration
     \[     0=\dSlam[0,\blam]\subset \dSlam[1,\blam]\subset
                        \dots\subset \dSlam[y-1,\blam]\subset \dSlam[y, \blam]
                         =\iRes\dSlam
     \]
    such that
    $\dSlam[r,\blam]/\dSlam[r-1,\blam]\cong q^{d_\blam^{B_r}}\dSlam[\blam{-}B_r]$
    as $\Rn[\beta]$-modules for $1\le r\le y$.
  \end{Corollary}

  \begin{proof}
    Rather than repeating the arguments of this
    section using the $\psi'$-basis, which we have not defined, we
    deduce the result from \autoref{C:GradedRestriction}. First, if~$M$
    is an $\Rn$-module let $M^\circledast=\Hom_{\AA}(M,\AA)$ be the
    graded dual of~$M$, where $\AA$ is in degree~$0$. Define the
    \textbf{defect} of~$\blam$ to be the non-negative integer
    $\defect\blam = (\Lambda, \beta_\blam)-\tfrac12(\beta_\blam,\beta_\blam)$,
    where $\bi^\blam=(i^\blam_1,\dots,i^\blam_n)$ and
    $\beta_\blam = \sum_{r=1}^n\alpha_{i^\blam_r}\in Q^+$. We can now prove the
    corollary.

    Let $\aset{\Slam[r,\blam]|0\le r\le y}$ be the filtration of $\Slam$
    given by \autoref{C:GradedRestriction} and set
    \[\dSlam[r,\blam]=q^{\defect\blam}(\Slam/\Slam[y-r,\blam])^\circledast,
               \qquad\text{for }0\le r\le y.\]
    Then $0=\dSlam[0,\blam]\subset \dSlam[1,\blam]\subset
                       \dots\subset \dSlam[y-1,\blam]\subset \dSlam[y, \blam]
                        =\iRes\dSlam
    $
    is a filtration of~$\dSlam$ as an $\Rn[n-1]$-module. Fix $r$ with $1\le r\le y$.
    By \autoref{C:GradedRestriction},
    $\Slam[y-r+1,\blam]/\Slam[y-r]\cong q^{d^\blam_{B_r}}\Slam[\blam-B_r]$
    as $\Rn[n-1]$-modules, so there is a short exact sequence
    \[    0\longrightarrow q^{d^\blam_{B_r}}\Slam[\blam-B_r]
           \longrightarrow \Slam/\Slam[y-r,\blam]
           \longrightarrow \Slam/\Slam[y-r+1,\blam]
           \longrightarrow0.
    \]
    By \cite[Theorem~7.25]{KMR:UniversalSpecht}, if $\bmu\in\Parts[m]$
    then $\dSlam[\bmu]\cong q^{\defect\bmu}(\Slam[\bmu])^\circledast$ as $\Rn[m]$-modules.
    Therefore, taking duals and shifting the degree by $q^{\defect\blam}$, there is a
    short exact sequence
    \[    0\longrightarrow \dSlam[r-1,\blam]
           \longrightarrow \dSlam[r,\blam]
      \longrightarrow q^{\defect\blam-\defect(\blam-B_r)-d^\blam_{B_r}}\dSlam[\blam-B_r]
           \longrightarrow0.
    \]
    This completes the proof because
    $d^{B_r}_\blam=\defect\blam-\defect(\blam-B_r)-d^\blam_{B_r}$ by
    \cite[Lemma~3.11]{BKW:GradedSpecht}.
  \end{proof}

\subsection*{Acknowledgements}
I thank Fred Goodman for alerting me to the gap in the published proofs
of the restriction formula and Liron Speyer and the two referees for
their many comments and corrections. This research was partially supported by
the Australian Research Council.


\bibliographystyle{abbrv}

\end{document}